\documentclass{amsart}
\usepackage{filecontents}
\usepackage[ocgcolorlinks,linktoc=all,linkcolor=blue]{hyperref}
\hypersetup{citecolor=blue,linkcolor=red}
\usepackage{cleveref}
\crefname{lemma}{Lemma}{Lemmata}
\crefname{equation}{equation}{equations}
 
\newtheorem{theorem}{Theorem}
\newtheorem{lemma}[theorem]{Lemma}
\newtheorem{proposition}[theorem]{Proposition}

\theoremstyle{definition}

\theoremstyle{remark}
\newtheorem{remark}[theorem]{Remark}

\numberwithin{equation}{section}

\renewcommand{\~}{\tilde}

\newcommand{\sub}{\subset}

\renewcommand{\S}{\mathbb{S}}

\newcommand{\8}{\infty}

\renewcommand{\a}{\alpha}
\renewcommand{\b}{\beta}
\newcommand{\g}{\gamma}
\renewcommand{\d}{\delta}

\renewcommand{\k}{\kappa}

\renewcommand{\t}{\theta}

\newcommand{\p}{F}

\newcommand{\G}{\Gamma}

\newcommand{\del}{\partial}

\newcommand{\fr}[2]{\frac{#1}{#2}}

\newcommand{\pf}[1]{\begin{proof} #1 \end{proof}}
\newcommand{\eq}[1]{\begin{equation}\begin{alignedat}{2} #1 \end{alignedat}\end{equation}}

\newcommand{\br}[1]{\left(#1\right)}

\newcommand{\mc}{\mathcal}

\newcommand{\hp}{\hphantom}

\protected\def\ignorethis#1\endignorethis{}
\let\endignorethis\relax

\begin{document}
\title[Differential Harnack inequalities on the sphere]{Harnack inequalities for evolving hypersurfaces on the sphere}

\author[P. Bryan]{Paul Bryan}
\address{Mathematics Institute, University of Warwick
Coventry, CV4 7AL, England}
\email{p.bryan@warwick.ac.uk}
\author[M.N. Ivaki]{Mohammad N. Ivaki}
\address{Institut f\"{u}r Diskrete Mathematik und Geometrie, Technische Universit\"{a}t Wien,
Wiedner Hauptstr. 8--10, 1040 Wien, Austria}
\email{mohammad.ivaki@tuwien.ac.at}
\author[J. Scheuer]{Julian Scheuer }
\address{Albert-Ludwigs-Universit\"{a}t,
Mathematisches Institut, Eckerstr. 1, 79104
Freiburg, Germany}
\email{julian.scheuer@math.uni-freiburg.de}
\date{\today}

\dedicatory{}
\subjclass[2010]{}
\keywords{Fully nonlinear curvature flows, Harnack estimates}

\begin{abstract}
We prove Harnack inequalities for hypersurfaces flowing on the unit sphere by $p$-powers of a strictly monotone, 1-homogeneous, convex, curvature function $f$, $0<p\leq 1.$ If $f$ is the mean curvature, we obtain stronger Harnack inequalities.
\end{abstract}

\maketitle

\section{Introduction}
We consider the evolution of a family of embeddings
\[x:M^n\times[0,T)\to M_c\]
of a smooth, closed manifold $M^n$ by
\eq{\label{eq:CurvFlow}
\partial_tx = -\p\nu,
}
where \(M_c\) is the simply connected space form of constant sectional curvature \(c\geq 0\) and $\p\in C^{\8}(\G_+)$ is a strictly monotone,  symmetric function of the eigenvalues of the Weingarten map \(\mathcal{W}\) (principal curvatures) \(\kappa_1, \cdots, \kappa_n\). Strict monotonicity ensures the flow is parabolic. We will need to make some further assumptions on the speed to obtain Harnack inequalities, namely that \[\p=f^p,\quad 0<p\leq 1,\] where \(f\) is $1$-homogeneous and convex. 

Under these assumptions, our principal results are Harnack inequalities for flows of strictly convex hypersurfaces on the sphere. These results extend the Harnack inequalities obtained in \cite{BryanIvaki:08/2015,BryanLouie:04/2016} on the sphere to a broader class of flows, similar to, though more restrictive than the class of flows in Euclidean space \cite{Andrews:09/1994} for which Harnack inequalities are known. We obtain the following theorem.

\begin{theorem}\label{thm:harnack}
Let $f$ be a strictly monotone, $1$-homogeneous and convex curvature function, $0<p\leq 1$, and let $\p=f^p$. Let $x$ be a solution to \eqref{eq:CurvFlow} in a simply connected space form of constant sectional curvature \(c\geq 0\), and such that $M_t=x(M,t)$ is strictly convex for all $t.$ Then $F$ satisfies
\[
\partial_t \p-b^{ij}\nabla_i\p\nabla_j\p+\frac{p\p}{(p+1)t}>0.
\]
For $f=H$ the following stronger estimates hold:
If $\frac{1}{2}+\frac{1}{2n}< {p}< 1,$ then
\[
\partial_t H^{p} - b^{ij}\nabla_iH^{p}\nabla_jH^{p} - \frac{c {p}}{2{p}-1}H^{2{p}-1} + \frac{{p}}{{p}+1} \frac{H^{p}}{t} > 0
\]
and if $0<{p}\leq \frac{1}{2} + \frac{1}{2n}$ or $p=1$, then
\[
\partial_t H^{p} - b^{ij}\nabla_iH^{p}\nabla_jH^{p} - c n{p}H^{2{p}-1} + \frac{{p}}{{p}+1} \frac{H^{p}}{t} > 0.
\]
\end{theorem}

A number of authors have studied Harnack inequalities in Euclidean space. The genesis of such study is \cite{Hamilton:/1995} where Hamilton proves a Harnack inequality for the mean curvature flow of convex hypersurfaces. Harnack inequalities for other flows have been obtained in \cite{Chow:06/1991,Ivaki:11/2015,Ivaki:09/2015,Li:/2011,Smoczyk:/1997}, including flows by powers of Gauss curvature, centro-affine normal flows, and flows by powers of inverse mean curvature. The most general results, subsuming many other results, were obtained in \cite{Andrews:09/1994} for so-called \(\alpha\)-convex and \(\alpha\)-concave speeds. The philosophy as espoused in \cite{Hamilton:/1995} by Hamilton is that equality should be attained by self-similar solutions, a.k.a. solitons, as originally motivated by the fact that the equality case in the seminal work of Li-Yau \cite{LiYau:/1986} is attained on the heat kernel, a self-similar solution of the heat equation. Such self-similarity leads one to study the Harnack quadratic,
\begin{equation}
\label{eq:harnack_quadratic}
Q = \partial_t \p - b^{ij} \nabla_i \p \nabla_j \p,
\end{equation}
where $\{b^{ij}\}$ is the inverse of the second fundamental form of a strictly convex hypersurface. This quadratic also arises, seemingly magically when changing parametrization from the ``Gauss map parametrization'' (where the calculations are almost trivial), to the ``Standard parametrization'' \cite{Andrews:09/1994}. Several authors, beginning with \cite{ChowChu:06/2001} have also investigated this quadratic, relating it to the second fundamental form of a degenerate metric on space-time \cite{HelmensdorferTopping:01/2013,Kotschwar:/2009}.

In the sphere, we do not have quite the same notion of self-similarity and the Gauss map parametrization does not seem to have quite the same ``magical'' properties as in Euclidean space. Using the ``Euclidean'' Harnack quadratic \eqref{eq:harnack_quadratic} on the sphere, we immediately encounter new difficulties arising from the background curvature introducing the ``remainder term'' $R$ of Proposition \ref{thm:Evchi}. In addition to the positivity required for Euclidean Harnack inequalities, on the sphere we also require positivity of $R$; therefore, our Harnack inequalities apply to a restricted class of flows as compared with \cite{Andrews:09/1994}. In the particular case when the speed is a function of the mean curvature, \(H^p\), the computation becomes tractable (\Cref{RSphere}, \Cref{lemma : lem9}), and suitably modifying the Harnack quadratic, we can cancel some bad terms to obtain the second Harnack inequality in \Cref{thm:harnack}.

Let us remark in passing that our computations recover most of the Harnack inequalities in Euclidean space mentioned above. In the case of space-forms of negative curvature, we find that essentially all the terms have the wrong sign and no Harnack inequality seems possible.


This paper is laid out as follows: in \Cref{sec:prelim} we define our notational conventions and recall some standard definitions and identities. In \Cref{sec:basic_evolution} we give some standard evolution equations and commutators and carry out the tedious task of computing the evolution of various quantities necessary for the main argument. \Cref{sec:main_evolution} combines these computations into evolution equations for the Harnack quadratics we study. Then, applying these calculations, we derive the Harnack inequalities in \Cref{sec:Harnack}. In this section, we present several variants depending on the strength of our assumptions. To finish, we prove preservation of convexity in \Cref{sec:convexity} for various flows in order to show that the assumption of convexity is reasonable.

\section*{Acknowledgment}

The authors would like to thank Knut Smoczyk and the Institut f\"{u}r Differentialgeometrie at Leibniz Universit\"at for hosting a research visit where part of this work took place. Bennet Chow was also very encouraging, suggesting that Harnack inequalities appear to be quite robust and should hold for a broad class of curvature flows, inspiring us to undertake the involved calculations required. The first author would also like to thank the third author for the gift of a bottle of French wine, only to be opened upon completion of this paper which served as strong motivation for completion. The first author was a Riemann Fellow at the Riemann Center for Geometry and Physics, Leibniz Universit\"{a}t whilst this work was conducted. The work of the second author was supported by Austrian Science Fund (FWF) Project M1716-N25 and the European Research Council (ERC) Project 306445. 

\section{Preliminaries}
\label{sec:prelim}

For a general Riemannian manifold $(M^n,g)$ let $\nabla$ be the Levi-Civita connection for $g$. Let $(\partial_i)$, $1\leq i\leq n$, be a coordinate frame.
We shall write \(\nabla_i=\nabla_{\partial_i}\) for covariant derivatives in direction $\partial_i$ and also use the notation \(\nabla^i = g^{ik} \nabla_k\). The Christoffel symbols are defined by 
\[\nabla_{\partial_j}\partial_i=\G^k_{ij}\partial_k.\]
For a $(k,l)$ tensor, $(\nabla_i T)^{i_1\cdots i_k}_{j_1\cdots j_l}$ will be written $\nabla_i T^{i_1\cdots i_k}_{j_1\cdots j_l}$. Second covariant derivatives will be written 
\[\nabla^2_{ij} = \nabla_{\partial_i} \nabla_{\partial_j} - \nabla_{\nabla_{\partial_i}\partial_j}\] and $(\nabla^2)^i_j = g^{ik} \nabla^2_{kj}$.

We use Hamilton's convention \cite[p.~258]{Hamilton:/1982} for the Riemannian curvature tensor, namely
\eq{\label{Rm}{{R_{ij}}^l}_k\partial_l&=\left(\frac{\partial}{\partial x^i}\G^l_{jk}-\frac{\partial}{\partial x^j}\G^l_{ik}+\G^l_{ir}\G^r_{jk}-\G^l_{jr}\G^r_{ik}\right)\partial_l\\
		&=\nabla_{i}\nabla_{j}\partial_k-\nabla_{j}\nabla_{i}\partial_k\\
        &=\nabla^2_{ij}\partial_k-\nabla^2_{ji}\partial_k,}
where the whole relation \eqref{Rm} is also known as the Ricci identity. It follows that for a function $f\in C^1(M)$ we have
\eq{\label{Rc}\nabla^3_{ijk}f-\nabla^3_{jik}f={R_{ijk}}^l\nabla_lf,}
where \(\nabla f = df\) is the covariant derivative of \(f\) and the subscript \(i\) refers to the \(i\)'th component of a one-form; see \cite[p.~258]{Hamilton:/1982}.

Let $\bar{g}$ and $\bar{R}$ denote, respectively, the metric and the curvature tensor of $M_c$. Now we specify to the situation where \(M_t := x(M^n,t)\) and \(g = x_t^{\ast} \bar{g}\) denotes the time dependent induced metric on \(M\) with $\nabla$ the corresponding time dependent Levi-Civita connection. Write $\nu$ for the outer unit normal to $M_t$, which gives rise to a frame \(\{\partial_0 = \nu, (x_t)_{\ast} \partial_1, \cdots, (x_t)_{\ast} \partial_n\}\) on \(M_c\) in a neighborhood of \(M_t\). Let Greek indices range from \(0\) to \(n\) and Latin indices range from \(1\) to \(n\).

The Riemann curvature tensor of \(M_c\) satisfies \(\bar{R}_{\a\beta\gamma\delta} = c(\bar{g}_{\a\gamma}\bar{g}_{\beta\delta} - \bar{g}_{\a\delta}\bar{g}_{\beta\gamma})\). We may write the metric $g = \{g_{ij}\}$, the second fundamental form $A = \{h_{ij}\}$, the Weingarten map $\mathcal{W} = \{h^i_j\} = \{g^{mi} h_{jm}\}$ and the Riemann curvature tensor $\{R_{ijkl}\}$ with respect to the given frame.

The mean curvature of $M^n$ is the trace of the Weingarten map (equivalently the trace of the second fundamental form with respect to $g$), $H = g^{ij}h_{ij} = h^i_i$. We also use the following standard notation
\[
(h^2)_i^j = g^{mj}g^{rs}h_{ir}h_{sm},
\]
\[
(h^2)_{ij} = g_{kj} (h^2)_i^k = h^k_i h_{kj},
\]
\[
|A|^2 = g^{ij}g^{kl}h_{ik}h_{lj} = h_{ij}h^{ij}.
\]
Here, $\{g^{ij}\}$ is the inverse matrix of $\{g_{ij}\}.$ For a strictly convex hypersurface, \(A\) is strictly positive-definite and hence has a strictly positive-definite inverse, which we denote by
\[
b = \{b^{ij}\}.
\]
The relations between $A$, $R$, and $\bar{R}$ are given by the Gau{\ss} and Codazzi equations:
\[
\begin{split}
R_{ijkl} &= \br{x^{*}\bar{R}}_{ijkl} + h_{ik}h_{jl} - h_{il}h_{jk} \\
&= c({g}_{ik}{g}_{jl} - {g}_{il}{g}_{jk}) + h_{ik}h_{jl} - h_{il}h_{jk},
\end{split}
\]
\[
\nabla_i h_{jk} = \nabla_k h_{ij},
\]
valid for space forms. 
We will need some notation for derivatives of the speed \(\p\). Let us write
\[
\p^{i}_{j} = \frac{\partial \p}{\partial h^{j}_{i}}
\]
for the first partial derivatives of \(\p\). We may also think of \(\p\) as a function of the metric and second fundamental form
\[
\p(g, h) = \p(g^{ik} h_{kj}).
\]
From this point of view, for the first and second partial derivatives, let us write
\[
\p^{ij} = \frac{\partial\p}{\partial h_{ij}}, \quad \p^{ij,kl} = \fr{\partial^2\p}{\partial h_{kl} \partial h_{ij}}.
\]
The trace of $\{\p^{ij}\}$ with respect to the metric will be written
\[
\operatorname{tr}(\dot{\p}) = g_{ij} \p^{ij}.
\]
Let us define the operator
$
\Box = \p^{ij} \nabla^2_{ij},
$
i.e., for a $(k,l)$-tensor $T$, $\Box T$ reads in coordinates 
\eq{\Box T_{j_1\dots j_l}^{i_1\dotsi_k}=F^{rs}(\nabla^2_{rs}T)^{i_1\dots i_k}_{j_1\dots j_l}.}

The $\Box$ operator satisfies the product rule. For smooth functions $\phi$ and $\psi$ we have
\begin{equation}
\label{eq:productbox}
\Box (\phi\psi) = \phi \Box \psi + \psi \Box \phi + 2 \p^{ij} \nabla_i \phi \nabla_j\psi.
\end{equation}

We frequently make use, without comment, of the formula for differentiating an inverse
\[
\frac{\partial g^{ij}}{\partial g_{kl}} = - g^{kj} g^{il}.
\]
First derivatives of \(\p\) from the two perspectives are related by
\begin{equation}
\label{eq:delh}
\p^{ij} = \frac{\partial \p}{\partial h_l^k} \frac{\partial h_l^k}{\partial h_{ij}} = \p^l_k g^{ik} \delta^j_l = g^{ik} \p^j_k
\end{equation}
and
\begin{equation}
\label{eq:delg}
\frac{\partial\p}{\partial g_{ij}} = \p^{l}_{k} \frac{\partial h^{k}_{l}}{\partial g_{ij}} = -\p^{l}_{k} g^{ki} g^{rj} h_{rl} = -\p^{li}h^{j}_{l}.
\end{equation}
We will also need the mixed second derivatives,
\begin{equation}
\label{eq:delhdelg}
\begin{split}
\frac{\partial \p^{ij}}{\partial g_{kl}} &= \frac{\partial}{\partial g_{kl}} \left(g^{sj} \p^{i}_{s} \right) = - g^{kj}g^{sl} \p^{i}_{s} - g^{sj} (\p^i_s)^{mk} h^l_m \\
&= - g^{kj} \p^{il} - g^{sj} g_{ns} \p^{in,mk} h^l_m \\
&= - g^{kj} \p^{il} - \p^{ij,mk} h^l_m,
\end{split}
\end{equation}
where we applied \eqref{eq:delg} to \(\p^i_s\) in the first line.

Covariant derivatives of \(\p\) satisfy
\begin{equation}
\label{eq:delphi}
\nabla_k \p = \p^{ij} \nabla_k h_{ij}
\end{equation}
and the covariant derivative of the trace satisfies
\begin{equation}
\label{eq:delPhi}
\nabla_k \operatorname{tr}(\dot{\p}) = g_{ij} \p^{ij,rs} \nabla_k h_{rs}.
\end{equation}

\section{Basic Evolution equations}
\label{sec:basic_evolution}

The evolution equation derived in this section hold for a general curvature function $F$. In the whole paper we only consider flows of strictly convex hypersurfaces.

Following \cite{Andrews:09/1994, Chow:06/1991, Hamilton:/1995, Smoczyk:/1997}, in this section we collect basic evolution equations that are needed to calculate the evolution of the quantities
\[
\chi_1 =t(\partial_t \p- b^{ij} \nabla_i \p \nabla_j \p) +\delta\p
\]
and
\[
\chi_2 =t(\partial_t \p - b^{ij} \nabla_i \p \nabla_j \p - c \p \operatorname{tr}(\dot{\p})) +\delta\p,
\]
where \(\delta \ne 0\) is an arbitrary, non-zero constant. The evolution equation of $\chi_1$ will be used for obtaining Harnack estimates for powers of convex 1-homogeneous curvature speeds $\p$ and the evolution equation of $\chi_2$ will be employed for obtaining stronger Harnack estimates for flows by powers of the mean curvature $\p(H)=H^{p}$ with $p\in(0,1].$ Note that in Euclidean space $\chi_1=\chi_2.$

Let us make a few definitions to keep the calculations more manageable. Let
\[
\a_{ij} = \nabla^2_{ij} \p + \p(h^2)_{ij}, \quad \gamma_{ij} = b^{kl} \nabla_k \p \nabla_l h_{ij}, \quad \eta_{ij} = \a_{ij} - \gamma_{ij}
\]
and define
\[
\beta = \p^{ij} \a_{ij} = \Box\p + \p \p^{ij}(h^2)_{ij}, \quad \theta =  b^{ij} \nabla_i \p \nabla_j \p,
\]
so that from the evolution of \(\p\) below (\cref{lem:evolution}, \cref{eq:delt_speed}) we may write our main Harnack quantities as
\[
\chi_1 = t(\partial_t\p- \theta) + \delta\p
\quad \mbox{and}\quad
\chi_2 = t(\beta - \theta) + \delta\p.
\]
We begin by recalling some standard evolution equations and commutators and then break the calculation into several lemmas.

The evolution equations in the following lemma are standard and can be found in many places \cite{Andrews:09/1994, Chow:06/1991, Hamilton:/1995, Huisken:/1987a, Smoczyk:/1997}. The necessary tools are commuting derivatives, using the definition of the curvature tensor for space forms, the Gauss equation, and the Codazzi equation as described in the previous section. Compare also \cite[p.~94-95]{Gerhardt:/2006} and the formula \cite[eq.~(6.17)]{Gerhardt:01/1996}.
\begin{lemma}
\label{lem:evolution}
The following evolution equations hold
\begin{enumerate}
\item \label{eq:delt_metric} $\partial_tg_{ij} = -2\p h_{ij},$
\item \label{eq:delt_inversemetric} $\partial_t g^{ij} = 2\p h^{ij},$
\item \label{eq:delt_sff} $\partial_t h_{ij} = \nabla^2_{ij} \p - \p(h^2)_{ij} + c \p g_{ij},$
\item \label{eq:delt_weingarten} $\partial_t h_i^j = (\nabla^2)^j_i\p + \p(h^2)_i^j + c \p\delta_i^j = \a^j_i + c \p\delta_i^j,$
\item \label{eq:delt_sff_box} \begin{align*}
\partial_t h_{ij} &= \Box h_{ij} + \p^{kl} (h^2)_{kl} h_{ij} - (\p^{kl}h_{kl} + \p) (h^2)_{ij} \\
& \quad + \p^{kl,rs}\nabla_i h_{kl}\nabla_j h_{rs} \\
& \quad + c \{(\p + \p^{kl}h_{kl}) g_{ij} - \operatorname{tr}(\dot{\p}) h_{ij}\},
\end{align*}
\item \label{eq:delt_weingarten_box} \begin{align*}
\partial_t h_i^j &= \Box h_i^j + \p^{kl} (h^2)_{kl} h_i^j - (\p^{kl}h_{kl} - \p) (h^2)_i^j \\
& \quad + \p^{kl,rs}\nabla_i h_{kl}\nabla^j h_{rs} \\
& \quad + c \{(\p + \p^{kl}h_{kl}) \delta_i^j - \operatorname{tr}(\dot{\p}) h_i^j\},
\end{align*}
\item \label{eq:delt_inversesff} \begin{align*}
\partial_t b^{ij} &= \Box b^{ij} - \p^{rs} (h^2)_{rs} b^{ij} + (\p^{kl}h_{kl} + \p) g^{ij} \\
& \quad - \left(2b^{lq}\p^{kp} + \p^{kl,pq}\right) b^{ir}b^{js} \nabla_r h_{kl} \nabla_s h_{pq} \\
& \quad - c \{(\p + \p^{kl}h_{kl}) b^{ir}b^{j}_{r} - \operatorname{tr}(\dot{\p}) b^{ij}\},
\end{align*}
\item \label{eq:delt_squaredsff} $\partial_t (h^2)_{ij} = h^k_j \nabla^2_{ik} \p + h^k_i \nabla^2_{jk} \p  + 2c\p h_{ij},$
\item \label{eq:delt_speed} $\partial_t \p = \Box \p + \p\p^{ij}(h^2)_{ij} + c \p\p^{ij}g_{ij} = \beta + c\p\operatorname{tr}(\dot{\p}).$
\end{enumerate}
\end{lemma}

\begin{lemma}
\label{EvGamma}
The Christoffel symbols evolve according to
\begin{equation}
\partial_t {\G}^{k}_{ij} = -\p g^{kl} \nabla_l h_{ij} - g^{kl} h_{li} \nabla_j \p - g^{kl} h_{lj} \nabla_i \p + g^{kl} h_{ij} \nabla_l \p.
\end{equation}
\end{lemma}

\begin{proof}
In local coordinates, we have
\[
\G^{k}_{ij} = \frac{1}{2} g^{kl} \left(\partial_j g_{il} + \partial_i g_{jl} - \partial_l g_{ij}\right).
\]
Since $\partial_t \G^{k}_{ij}$ is a tensor, we may calculate using normal coordinates at any given point, at which \(\G^k_{ij} = 0\). Then we have
\[
\frac{1}{2} \partial_t g^{kl} (\partial_j g_{il} + \partial_i g_{jl} - \partial_l g_{ij}) = 2\p g^{kr} h_{rs} \G^s_{ij} = 0
\]
from \cref{lem:evolution}, \cref{eq:delt_inversemetric}. Now commuting derivatives \([\partial_t, \partial_i] = 0\), and using the Codazzi equations we obtain
\[
\begin{split}
\partial_t {\G}^{k}_{ij} & = \frac{1}{2} g^{kl} \left(\partial_j \partial_t g_{il} + \partial_i \partial_t g_{jl} - \partial_l \partial_t g_{ij}\right) \\
&= - g^{kl} \left(\partial_j (\p h_{il}) + \partial_i (\p h_{jl}) - \partial_l (\p h_{ij})\right) \\
&= -\p g^{kl} \partial_l h_{ij} - g^{kl} h_{il} \partial_j \p  - g^{kl} h_{lj} \partial_i \p + g^{kl} h_{ij} \partial_l \p.
\end{split}
\]
The result follows since in normal coordinates, \(\nabla_i = \partial_i\) at our given point.
\end{proof}

We require the commutators \([\nabla, \Box]\) and \([\partial_t, \Box]\). Without further comment we will also use the fact that \([\partial_t, \nabla] f = 0\) for any smooth function \(f\).

\begin{lemma}
\label{lem:gradBox}
For every smooth function $f$, the commutation relation holds
\[
\begin{split}
([\nabla, \Box]f)_i &= \nabla_i \Box f - (\Box \nabla f)_i\\
                &= \p^{kl,rs} \nabla_i h_{rs} \nabla^2_{kl} f + \p^{kl}\left(h^{m}_{k}h_{li} - h_{kl}h^{m}_{i}\right) \nabla_m f \\
&\quad + c \p^{kl} g_{li} \nabla_k f - c \operatorname{tr}(\dot{\p}) \nabla_i f
\end{split}.
\]

\end{lemma}
\begin{proof}
From the Ricci identity \eqref{Rc} we get

\[
\nabla^3_{ikl}f-\nabla^3_{kli}f={{R}_{ikl}}^m\nabla_m f
\]
and thus we obtain
\[
\nabla_i (\p^{kl} \nabla^2_{kl} f) - (\p^{kl}(\nabla^2_{kl} \nabla f))_i = \p^{kl,rs} \nabla_i h_{rs} \nabla^2_{kl}f+ \p^{kl}{{R}_{ikl}}^{m} \nabla_m f.
\]
From the Gauss equation we obtain
\[
\begin{split}
{{R}_{ikl}}^{m} \nabla_m f &= \left(c\left(g^{pm}g_{il}g_{kp}  - g^{pm}g_{pi}g_{kl}\right) + g^{pm} h_{il}h_{kp} - g^{pm}h_{kl}h_{pi}\right) \nabla_m f \\
&= c\left(g_{li} \nabla_k f - g_{kl} \nabla_i f\right) + \left(h^{m}_{k}h_{li} - h_{kl}h^{m}_{i}\right) \nabla_m f.
\end{split}
\]
\end{proof}
\begin{lemma}
\label{lem:deltBox}
The following commutation relation holds
\[
\begin{split}
[\partial_t, \Box] \p &= (\partial_{t}\Box - \Box\partial_{t}) \p = \p^{ij,kl} \nabla^2_{ij} \p (\a_{kl} + c \p g_{kl}) \\
&\quad + 2\p^{ij}h^{k}_{i} (\p \nabla^2_{k
j} \p + \nabla_k \p \nabla_j \p) + (\p - \p^{ij}h_{ij})| \nabla\p|^{2}.
\end{split}
\]
\end{lemma}
\begin{proof}
First, let us calculate the evolution of \(\p^{ij}\), which will also prove useful later. From the mixed derivative \cref{eq:delhdelg}, the evolution of the metric (\cref{lem:evolution}, \cref{eq:delt_metric}), and the evolution of the second fundamental form (\cref{lem:evolution}, \cref{eq:delt_sff}) we compute
\begin{equation}
\label{eq:deltBox}
\begin{split}
\partial_{t} \p^{ij} &= \p^{ij,kl} \partial_t h_{kl} + \frac{\partial\p^{ij}}{\partial g_{kl}} \partial_t g_{kl} \\
&= \p^{ij,kl} \left(\nabla^2_{kl} \p - \p(h^2)_{kl} + c \p g_{kl}\right) + 2\p \p^{ij,kl} h_{lm}h^{m}_{k} + 2\p\p^{jk}g^{li}h_{kl} \\
&= \p^{ij,kl} \left(\nabla^2_{kl} \p + \p(h^2)_{kl} + c \p g_{kl}\right) + 2\p\p^{jk}h^{i}_{k} \\
&= \p^{ij,kl} \left(\a_{kl} + c \p g_{kl}\right) + 2\p\p^{jk}h^{i}_{k}.
\end{split}
\end{equation}
Next, the commutator of \(\partial_t\) and \(\nabla^2_{ij}\) is given by,
\begin{equation}
\label{eq:deltnabla2}
\begin{split}
\left(\partial_{t}\nabla^2_{ij} - \nabla^2_{ij}\partial_{t}\right) \p &= \partial_t \left(\nabla_i \nabla_j \p - \nabla_{\nabla_i \partial_j} \p\right) - \nabla_i \nabla_j \partial_t \p + \nabla_{\nabla_i \partial_j} \partial_t \p \\
&= - \partial_t \left(\G_{ij}^k \nabla_{\partial_k} \p\right) + \G_{ij}^k \nabla_{\partial_k} \partial_t \p \\
&= - \nabla_k \p \partial_t \G_{ij}^k.
\end{split}
\end{equation}
We obtain from \eqref{eq:deltBox}, \eqref{eq:deltnabla2}, and the evolution of the Christoffel symbols \cref{EvGamma},
\[
\begin{split}
\left(\partial_{t}\Box - \Box\partial_{t}\right) \p &= \left(\partial_{t}\p^{ij}\right) \nabla^2_{ij} \p - \p^{ij}\nabla_k \p\partial_t \G^{k}_{ij} \\
&= \left(\p^{ij,kl} \left(\a_{kl} + c \p g_{kl}\right) + 2\p\p^{jk}h^{i}_{k}\right) \nabla^2_{ij} \p \\
&\quad + \p^{ij} \nabla_k \p \left(\p g^{kl} \nabla_l h_{ij} + h^k_i \nabla_j \p + h^k_j \nabla_i \p - g^{kl} h_{ij} \nabla_l \p\right) \\
&= \p^{ij,kl} \nabla^2_{ij} \p \left(\a_{kl} + c \p g_{kl}\right) \\
&\quad + 2\p^{jk}h^{i}_{k}\p \nabla^2_{ij} \p + h^k_i \p^{ij} \nabla_k \p \nabla_j \p + h^k_j \p^{ij} \nabla_k \p \nabla_i \p \\
&\quad + \p g^{kl}\nabla_k \p \p^{ij} \nabla_l h_{ij} - \p^{ij} h_{ij} g^{kl}\nabla_k \p \nabla_l \p.
\end{split}
\]
The result now follows from \(|\nabla \p|^2 = g^{kl}\nabla_k \p \nabla_l \p\) and \(\nabla_l \p = \p^{ij} \nabla_l h_{ij}\) by the chain rule.
\end{proof}
The next ingredient is the evolution of the covariant derivative, \(\nabla \p = d\p\).
\begin{lemma}
\label{lem:Evgradphi}
There holds
\[
\begin{split}
\left((\partial_{t}-\Box)\nabla\p\right)_{i} &= \p^{kl,rs}\nabla_i h_{rs} \a_{kl} + 2 \p^{kl} b^{rs} \p(h^2)_{rl} \nabla_i h_{ks} \\
&\quad + \p^{kl}(h^2)_{kl}\nabla_i \p + \left(\p^{kl}h^{m}_{l}h_{ki} - \p^{kl}h_{kl}h^{m}_{i}\right) \nabla_m \p\\
&\quad + c\left(\p^{kl}g_{ki} \nabla_l \p + \p \nabla_i \operatorname{tr}(\dot{\p})\right).
\end{split}
\]
\end{lemma}
\begin{proof}
Using the evolution of \(\p\) (\cref{lem:evolution}, \cref{eq:delt_speed}) and the commutator \([\nabla, \Box]\) from \cref{lem:gradBox}, we compute
\[
\begin{split}
\partial_{t}\nabla_i \p - (\Box\nabla \p)_{i} &= \nabla_i \partial_t \p - \nabla_i \Box \p + ([\nabla, \Box] \p)_i \\
&= \nabla_i \left(\Box\p + \p^{kl}(h^2)_{kl}\p + c \operatorname{tr}(\dot{\p})\p\right) - \nabla_i (\Box\p) \\
&\quad + \p^{kl,rs} \nabla_i h_{rs} \nabla^2_{kl} \p + c\p^{kl}g_{ki} \nabla_l \p - c\operatorname{tr}(\dot{\p})\nabla_i \p \\
&\quad + (\p^{kl}h^{m}_{l}h_{ki} - \p^{kl}h_{kl}h^{m}_{i}) \nabla_m \p \\
&= \p^{kl}(h^2)_{kl}\nabla_i \p + \p^{kl,rs}\nabla_i h_{rs} (h^2)_{kl}\p \\
&\quad+ \p\p^{kl}(h^s_l \nabla_i h_{ks} + h^r_k \nabla_i h_{rl})  + c \p\nabla_i\operatorname{tr}(\dot{\p})\\
&\quad+ \p^{kl,rs} \nabla_i h_{rs} \nabla^2_{kl} \p + c\p^{kl}g_{ki}\nabla_l \p \\
&\quad + (\p^{kl}h^{m}_{l}h_{ki} - \p^{kl}h_{kl}h^{m}_{i}) \nabla_m \p \\
&= \p^{kl,rs}\nabla_i h_{rs} \left((h^2)_{kl}\p + \nabla^2_{kl} \p\right) + 2 \p\p^{kl} h^s_l \nabla_i h_{ks} \\
&\quad + \p^{kl}(h^2)_{kl}\nabla_i \p + (\p^{kl}h^{m}_{l}h_{ki} - \p^{kl}h_{kl}h^{m}_{i}) \nabla_m \p \\
&\quad + c \left(\p\nabla_i\operatorname{tr}(\dot{\p}) + \p^{kl}g_{ki}\nabla_l \p\right),
\end{split}
\]
where in the third equality we used
$
\nabla_i (h^2)_{kl} = \nabla_i (g^{sr} h_{ks} h_{rl}) = h^s_l \nabla_i h_{ks} + h^r_k \nabla_i h_{rl}.
$
We obtain the result, since 
$
b^{rs} (h^2)_{rl} = b^{rs} h_{rm} h^m_l = \delta^s_m h^m_l = h^s_l.
$
\end{proof}
Now we may proceed to the calculations of \(\partial_t \beta\) and \(\partial_t \theta\).
\begin{lemma}
\label{lem:evbeta}
The quantity $\beta$
satisfies
\[
\begin{split}
(\partial_{t} - \Box)\beta &= \left(\p^{ij}(h^2)_{ij} + c\operatorname{tr}(\dot{\p}) \right)\beta + (\p - \p^{ij}h_{ij}) |\nabla\p|^{2}+ \p^{ij,kl} \a_{ij} \a_{kl}\\ 
&\quad + 2\p^{ij}h^{k}_{i}\nabla_k \p \nabla_j \p + 2b^{il}\p^{jk} (2\nabla^2_{ij}\p\p(h^2)_{kl} + \p(h^2)_{ij}\p(h^2)_{kl}) \\
&\quad + cR_{\beta},
\end{split}
\]
where
$
R_{\beta} = \p \Box \operatorname{tr}(\dot{\p}) + 2\p^{kl} \nabla_k \operatorname{tr}(\dot{\p}) \nabla_l \p + \p \p^{ij,kl}g_{kl} \a_{ij} + 2\p^{2}\p^{ij}h_{ij}.
$
\end{lemma}
\begin{proof}
Let us break up the calculation of
\[
(\partial_{t} - \Box)\beta =  \partial_{t}\Box\p + \partial_{t} (\p\p^{ij} (h^2)_{ij}) - \Box\Box\p - \Box(\p\p^{ij} (h^2)_{ij})
\]
into smaller pieces. First, we have lots of nice cancellation. Using the evolution of \(\p\) from \cref{lem:evolution}, \cref{eq:delt_speed} and the commutator relation from \cref{lem:deltBox} we have,
\begin{equation}
\label{eq:deltbeta1}
\begin{split}
&\partial_{t}\Box\p - \Box\Box\p - \Box(\p\p^{ij} (h^2)_{ij})\\ =\ &\Box\partial_t\p - \Box\Box\p - \Box(\p\p^{ij} (h^2)_{ij}) + [\partial_t, \Box] \p \\
=\ &\Box(\Box \p + \p\p^{ij}(h^2)_{ij} + c \p\operatorname{tr}(\dot{\p})) - \Box\Box\p - \Box(\p\p^{ij} (h^2)_{ij})\\
\quad &+ \p^{ij,kl} \nabla^2_{ij} \p (\a_{kl} + c \p g_{kl}) + 2\p^{ij}h^{k}_{i} (\p \nabla^2_{kj} \p + \nabla_k \p \nabla_j \p)\\
\quad &+ (\p - \p^{ij}h_{ij})| \nabla\p|^{2} \\
=\ &c (\operatorname{tr}(\dot{\p}) \Box \p + \p \Box \operatorname{tr}(\dot{\p}) + 2 \p^{kl} \nabla_k \operatorname{tr}(\dot{\p}) \nabla_l \p) \\
\quad &+ \p^{ij,kl} \nabla^2_{ij} \p (\a_{kl} + c \p g_{kl}) + 2\p^{ij}b^{kl} \p (h^2)_{il} \nabla^2_{kj} \p\\
\quad  &+ 2\p^{ij}h^{k}_{i} \nabla_k \p \nabla_j \p + (\p - \p^{ij}h_{ij})| \nabla\p|^{2} \\
=\ &c \operatorname{tr}(\dot{\p}) \Box \p + (\p - \p^{ij}h_{ij})| \nabla\p|^{2}  + \p^{ij,kl} \nabla^2_{ij} \p \a_{kl} \\
\quad &+ 2\p^{ij}b^{kl} \p (h^2)_{il} \nabla^2_{kj} \p + 2\p^{ij}h^{k}_{i} \nabla_k \p \nabla_j \p \\
\quad &+ c(\p \Box \operatorname{tr}(\dot{\p}) + 2 \p^{kl} \nabla_k \operatorname{tr}(\dot{\p}) \nabla_l \p + \p \p^{ij,kl} g_{kl} \nabla^2_{ij} \p)
\end{split}
\end{equation}
using, in the third equality, the product rule \eqref{eq:productbox} for \(\Box\)  and \(h^k_i = h^m_i b^{kl}h_{ml} = b^{kl} (h^2)_{il}\) since \(b\) is the inverse of \(A\).

Next from \eqref{eq:deltBox} and \cref{lem:evolution}, \cref{eq:delt_squaredsff} we obtain
\[
\begin{split}
\partial_{t} (\p^{ij}(h^2)_{ij}) &= \partial_{t}(\p^{ij}) (h^2)_{ij} + \p^{ij} \partial_t (h^2)_{ij} \\
&= \left(\p^{ij,kl} \left(\a_{kl} + c \p g_{kl}\right) + 2\p\p^{jk}h^{i}_{k}\right) (h^2)_{ij} \\
&\quad + \p^{ij} \left(h^k_j \nabla^2_{i
k} \p + h^k_i \nabla^2_{jk} \p + 2c\p h_{ij}\right) \\
&= \p^{ij,kl} (h^2)_{ij} \left(\a_{kl} + c \p g_{kl}\right) + 2\p\p^{jk} b^{il} (h^2)_{lk} (h^2)_{ij} \\
&\quad + 2 \p^{ij} b^{kl} (h^2)_{il} \nabla^2_{jk} \p  + 2c\p\p^{ij}h_{ij},
\end{split}
\]
again using \(h^k_i = b^{kl} (h^2)_{il}\) in the last equality.

The remaining term we need to compute is thus
\begin{equation}
\label{eq:deltbeta2}
\begin{split}
\partial_{t} (\p \p^{ij}(h^2)_{ij}) &= (\partial_{t} \p) \p^{ij}(h^2)_{ij} + \p \partial_t (\p^{ij} (h^2)_{ij}) \\
&= (\beta + c \p\operatorname{tr}(\dot{\p})) \p^{ij}(h^2)_{ij} \\
&\quad + \p \left(\p^{ij,kl} (h^2)_{ij} \left(\a_{kl} + c \p g_{kl}\right) + 2\p\p^{jk} b^{il} (h^2)_{lk} (h^2)_{ij} \right.\\
&\quad \left. + 2 \p^{ij} b^{kl} (h^2)_{il} \nabla^2_{jk} \p  + 2c\p\p^{ij}h_{ij}\right). \\
&= (\beta + c \p\operatorname{tr}(\dot{\p})) \p^{ij}(h^2)_{ij} + \p^{ij,kl} \p (h^2)_{ij} \a_{kl}  \\
&\quad + 2 \p\p^{ij} b^{kl} (h^2)_{il} \nabla^2_{jk} \p + 2\p^{jk} b^{il} \p (h^2)_{lk} \p (h^2)_{ij} \\
&\quad + c \left(\p \p^{ij,kl} g_{kl} \p (h^2)_{ij} + 2\p^2\p^{ij}h_{ij}\right).
\end{split}
\end{equation}

Now we add \eqref{eq:deltbeta1} and \eqref{eq:deltbeta2} together line by line to complete the proof.
\end{proof}

\begin{lemma}
\label{lem:Evtheta}
The quantity $\theta$ evolves according to
\[
\begin{split}
(\partial_{t} - \Box)\theta &= (\p^{ij}(h^2)_{ij} + c\operatorname{tr}(\dot{\p}))\theta \\
&\quad + (\p - \p^{ij}h_{ij})|\nabla\p|^{2} + 2\p^{ij}h^{k}_{i}\nabla_k\p\nabla_j\p \\
&\quad - \p^{kl,ij} (\gamma_{ij}\gamma_{kl}  - 2\a_{ij} \gamma_{kl}) - 2b^{il} \p^{jk} \left(\gamma_{ij} \gamma_{kl} - 2\a_{ij} \gamma_{kl} + \nabla^2_{ij}\p\nabla^2_{kl}\p\right) \\
&\quad + cR_{\theta},
\end{split}
\]
\end{lemma}
where
$
R_{\theta} = -(\p^{kl}h_{kl} + \p)b^{ir}b^{j}_{r}\nabla_i \p\nabla_j\p + 2 b^{j}_{k}\p^{kl}\nabla_l\p\nabla_j\p + 2 \p\p^{ij,kl} g_{ij} \gamma_{kl}.
$
\begin{proof}
Again using the product rule for \(\Box\), \eqref{eq:productbox}, we have
\begin{equation}
\label{eq:delt_theta}
\begin{split}
(\partial_{t} - \Box)\theta &= (\partial_{t}b^{ij} - \Box b^{ij})\nabla_i \p\nabla_j\p + b^{ij} ((\partial_{t} - \Box) (\nabla\p \otimes \nabla\p))_{ij} \\
&\quad - 2 \p^{kl} \nabla_k b^{ij} (\nabla_l (\nabla \p \otimes \nabla\p))_{ij} \\
&= (\partial_{t}b^{ij} - \Box b^{ij})\nabla_i \p\nabla_j\p + 2 b^{ij} ((\partial_{t} - \Box) (\nabla\p))_i \nabla_j\p\\
&\quad - 2 b^{ij} \p^{kl} \nabla^2_{ik} \p \nabla^2_{lj} \p - 4 \p^{kl} \nabla_k b^{ij} \nabla^2_{il} \p \nabla_j\p \\
&= (\partial_{t}b^{ij} - \Box b^{ij})\nabla_i \p\nabla_j\p + 2 b^{ij} ((\partial_{t} - \Box) (\nabla\p))_i \nabla_j\p \\
&\quad - 2 b^{ij} \p^{kl} \nabla^2_{ik} \p \nabla^2_{lj} \p + 4 \p^{kl} b^{ip}b^{jq} \nabla_k h_{pq} \nabla^2_{il} \p \nabla_j\p \\
&= (\partial_{t}b^{ij} - \Box b^{ij})\nabla_i \p\nabla_j\p + 2 b^{ij} ((\partial_{t} - \Box) (\nabla\p))_i \nabla_j\p \\
&\quad - 2 b^{ij} \p^{kl} \nabla^2_{ik} \p \nabla^2_{lj} \p + 4 \p^{kl} b^{ip}\gamma_{pk} \nabla^2_{il} \p,
\end{split}
\end{equation}
where in the second to last equality we used the formula for the derivative of the inverse \(b^{ij}\) of \(h_{ij}\) and the Codazzi equation in the last line, producing \(b^{jq} \nabla_k h_{pq} \nabla_j \p = b^{jq} \nabla_q h_{pk} \nabla_j \p = \gamma_{pk}\). The first term in the final line appears on the second to last line of the statement of the lemma (with indices relabelled). The second term is part of \(4 b^{il}\p^{jk} \a_{ij} \gamma_{kl}\) in the second to last line. So we must deal with the first two terms and show they add to the remainder of the statement. For the first term, we use the evolution of \(b^{ij}\) from \cref{lem:evolution}, \cref{eq:delt_inversesff} to calculate
\begin{equation}
\label{eq:delt_theta1}
\begin{split}
(\partial_{t}b^{ij} - \Box b^{ij})\nabla_i \p\nabla_j\p &= \nabla_i \p \nabla_j \p\left(-\p^{rs} (h^2)_{rs} b^{ij} + (\p^{kl}h_{kl} + \p) g^{ij} \right. \\
& \quad - \left(2b^{lq}\p^{kp} + \p^{kl,pq}\right) b^{ir}b^{js} \nabla_r h_{kl} \nabla_s h_{pq} \\
& \quad - \left. c \{(\p + \p^{kl}h_{kl}) b^{ir}b^{j}_{r} - \operatorname{tr}(\dot{\p}) b^{ij}\}\right) \\
&= \left(c\operatorname{tr}(\dot{\p}) - \p^{rs}(h)^2_{rs}\right) \theta+ (\p^{kl}h_{kl} + \p)|\nabla\p|^{2}\\ & \quad - (2b^{lq}\p^{kp} + \p^{kl,pq}) b^{ir}b^{js}\nabla_i\p\nabla_j\p\nabla_rh_{kl}\nabla_s h_{pq} \\
&\quad  - c(\p^{kl}h_{kl} + \p)b^{ir}b^{j}_{r}\nabla_i \p\nabla_j\p \\
&= \left(c\operatorname{tr}(\dot{\p}) - \p^{rs}(h)^2_{rs}\right) \theta  + (\p^{kl}h_{kl} + \p)|\nabla\p|^{2} \\
&\quad - \p^{kl,pq} \gamma_{kl} \gamma_{pq} - 2b^{lq}\p^{kp} \gamma_{kl} \gamma_{pq}  \\
&\quad -c(\p^{kl}h_{kl} + \p)b^{ir}b^{j}_{r}\nabla_i \p\nabla_j\p.
\end{split}
\end{equation}
For the second term, from the evolution of \(\nabla\p\) in \Cref{lem:Evgradphi}, we have
\begin{equation}
\label{eq:delt_theta2}
\begin{split}
&2 b^{ij} ((\partial_{t} - \Box) \nabla\p)_i \nabla_j\p\\
=\ &2 b^{ij} \nabla_j\p \big(\p^{kl,rs}\nabla_i h_{rs} \a_{kl} + 2 \p^{kl} b^{rs} \p(h^2)_{rl} \nabla_i h_{ks}  \\
\quad &+ \p^{kl}(h^2)_{kl}\nabla_i \p + \left(\p^{kl}h^{m}_{l}h_{ki} - \p^{kl}h_{kl}h^{m}_{i}\right) \nabla_m \p\\
\quad &+ c\left(\p^{kl}g_{ki} \nabla_l \p + \p \nabla_i \operatorname{tr}(\dot{\p})\right)\big) \\
=\ &2 b^{ij} \nabla_j\p \p^{kl}(h^2)_{kl}\nabla_i \p \\
\quad &- 2 b^{ij} \nabla_j\p \p^{kl}h_{kl}h^{m}_{i} \nabla_m \p + 2 b^{ij} \nabla_j\p \p^{kl}h^{m}_{l}h_{ki} \nabla_m \p \\
\quad &+ 2 b^{ij} \nabla_j\p \p^{kl,rs}\nabla_i h_{rs} \a_{kl}  + 4 b^{ij} \nabla_j\p \p^{kl} b^{rs} \p(h^2)_{rl} \nabla_i h_{ks} \\
\quad &+ c\left(2 b^{ij} \nabla_j\p \p^{kl}g_{ki} \nabla_l \p + 2 b^{ij} \nabla_j\p \p g_{kl}F^{kl,rs}\nabla_{i}h_{rs}\right) \\
=\ &2 \p^{kl}(h^2)_{kl}\theta - 2 \p^{kl}h_{kl} |\nabla\p|^2 + 2 \p^{kl} h^{m}_{l} \nabla_k\p \nabla_m \p \\
\quad &+ 2 \p^{kl,rs} \gamma_{rs} \a_{kl} + 4 b^{rs} \p^{kl} \gamma_{ks} \p(h^2)_{rl} \\
\quad &+ c\left(2 \p^{kl} b^j_k \nabla_j\p \nabla_l \p + 2 FF^{kl,rs}g_{kl}\gamma_{rs}\right),
\end{split}
\end{equation}
using the definitions of \(\theta, \a_{ij}\) and \(\gamma_{ij}\) as well as \(b^{ij}h_{ki} \nabla_j \p = \delta^j_k \nabla_j \p = \nabla_k \p\), and \(b^{ij} h^m_i = b^{ij} g^{mp}h_{pi} = \delta^j_p g^{mp} = g^{mj}\) in the last equality.

The proof is now completed by adding \eqref{eq:delt_theta1} and \eqref{eq:delt_theta2} line by line and adding also the final line from \eqref{eq:delt_theta}.
\end{proof}

\section{Main evolution equations}
\label{sec:main_evolution}

We start this section by calculating the evolution equations of $\chi_2$ and a slight modification, $\chi_3$, which will be employed to obtain Harnack estimates for flows by powers of the mean curvature. We will then focus on the evolution equation of $\chi_1$ which will enable us to deduce (weak) Harnack estimates for powers of 1-homogeneous convex speeds.
\begin{proposition}
\label{thm:Evchi}
Let $\delta \neq 0.$ For a general curvature function $F$ under flow \eqref{eq:CurvFlow} the quantity
$
\chi_2 = t(\beta - \theta) + \delta\p
$
satisfies
\eq{\label{thm:Evchi1}
\partial_t \chi_2 -\Box\chi_2 &= \left(\frac{\b-\t}{\delta\p} +\p^{ij}(h^2)_{ij} + c\operatorname{tr}(\dot{\p})\right)\chi_2 \\
& \quad + t\left(\p^{ij,kl} + 2b^{il}\p^{jk} - \frac{\p^{ij}\p^{kl}}{\delta\p}\right)\eta_{ij}\eta_{kl} + tc R,
}
where
\[
\eta_{ij} = \a_{ij} - \gamma_{ij} = \nabla^2_{ij}\p + (h^2)_{ij}\p - b^{rs}\nabla_r h_{ij}\nabla_s \p
\]
and
\begin{align}
\begin{split}\label{R}
R &= R_{\beta} - R_{\theta} \\
&= \p \Box \operatorname{tr}(\dot{\p}) + 2\p^{kl} \nabla_k \operatorname{tr}(\dot{\p}) \nabla_l \p + \p \p^{ij,kl}g_{kl} (\nabla^2_{ij} \p + \p(h^2)_{ij})\\
&\quad- 2\p \p^{ij,kl}g_{kl} b^{rs}\nabla_r h_{ij}\nabla_s \p + 2\p^{2}\p^{ij}h_{ij}\\
&\quad  +(\p^{kl}h_{kl} + \p)b^{ir}b^{j}_{r}\nabla_i \p\nabla_j\p - 2 b^{j}_{k}\p^{kl}\nabla_l\p\nabla_j\p.
\end{split}
\end{align}
\end{proposition}
\begin{proof}
We have
$
(\partial_t - \Box)\chi_2 = \beta - \theta + t(\partial_{t} - \Box)(\beta - \theta) + \delta(\partial_t \p - \Box\p).
$
First of all, the evolution equation for \(\p\), \cref{lem:evolution}, \cref{eq:delt_speed} gives us
\[
\delta(\partial_t \p - \Box\p) = \left(\p^{ij}(h^2)_{ij} + c\operatorname{tr}(\dot{\p})\right)\delta\p.
\]
Next, we note that
$
\p^{ij} \eta_{ij} = \beta - \theta
$
since \(\nabla_r \p = \p^{ij} \nabla_r h_{ij}\). Putting the two equations above together gives
\begin{equation}
\label{eq:deltchi1}
\begin{split}
&\beta - \theta + \delta(\partial_t \p - \Box\p) \\
=\ &\left(\frac{\beta-\theta}{\delta\p} + \p^{ij}(h^2)_{ij} + c\operatorname{tr}(\dot{\p})\right)\delta\p + t \frac{(\beta - \theta)^2}{\delta\p} - t \frac{(\p^{ij}\eta_{ij})^2}{\delta\p} \\
=\ &\frac{\beta-\theta}{\delta\p} \chi_2 + \left(\p^{ij}(h^2)_{ij} + c\operatorname{tr}(\dot{\p})\right)\delta\p  - t \frac{\p^{ij}\p^{kl}}{\delta\p} \eta_{ij}\eta_{kl}.
\end{split}
\end{equation}
The remaining term \(t(\partial_{t} - \Box)(\beta - \theta)\) is now just bookkeeping. Recall, \cref{lem:evbeta} states that
\begin{align*}
(\partial_{t} - \Box)\beta &= \left(\p^{ij}(h^2)_{ij} + c\operatorname{tr}(\dot{\p}) \right)\beta  & (A) \\
&\quad + (\p - \p^{ij}h_{ij}) |\nabla\p|^{2} + 2\p^{ij}h^{k}_{i}\nabla_k \p \nabla_j \p  & (B) \\
&\quad + \p^{ij,kl} \a_{ij} \a_{kl} & (C) \\
&\quad + 2b^{il}\p^{jk} (2\nabla^2_{ij}\p\p(h^2)_{kl} + \p(h^2)_{ij}\p(h^2)_{kl}) & (D) \\
&\quad + cR_{\beta}  & (E) \\
\intertext{while \cref{lem:Evtheta} states that}
(\partial_{t} - \Box)\theta &= (\p^{ij}(h^2)_{ij} + c\operatorname{tr}(\dot{\p}))\theta & (A') \\
&\quad + (\p - \p^{ij}h_{ij})|\nabla\p|^{2} + 2\p^{ij}h^{k}_{i}\nabla_k\p\nabla_j\p & (B') \\
&\quad - \p^{kl,ij} (\gamma_{ij}\gamma_{kl}  - 2\a_{ij} \gamma_{kl}) & (C') \\
&\quad - 2b^{il} \p^{jk} \left(\gamma_{ij} \gamma_{kl} - 2\a_{ij} \gamma_{kl} + \nabla^2_{ij}\p\nabla^2_{kl}\p\right) & (D') \\
&\quad + cR_{\theta}. & (E')
\end{align*}
Subtracting line by line, we have
\begin{align*}
(A) - (A') &= \left(\p^{ij}(h^2)_{ij} + c\operatorname{tr}(\dot{\p}) \right)\beta - \left(\p^{ij}(h^2)_{ij} + c\operatorname{tr}(\dot{\p}) \right)\theta \\
        &= \left(\p^{ij}(h^2)_{ij} + c\operatorname{tr}(\dot{\p}) \right)(\beta - \theta), \\
(B) - (B') &= (\p - \p^{ij}h_{ij}) |\nabla\p|^{2} + 2\p^{ij}h^{k}_{i}\nabla_k \p \nabla_j \p \\
    &\quad - (\p - \p^{ij}h_{ij})|\nabla\p|^{2} - 2\p^{ij}h^{k}_{i}\nabla_k\p\nabla_j\p \\
        &= 0, \\
(C) - (C') &= \p^{ij,kl} \a_{ij} \a_{kl} + \p^{kl,ij} (\gamma_{ij}\gamma_{kl}  - 2\a_{ij} \gamma_{kl})\\
    &= \p^{ij,kl} (\a_{ij} - \gamma_{ij}) (\a_{kl} - \gamma_{kl}) \\
&= \p^{ij,kl} \eta_{ij} \eta_{kl}, \\
(D) - (D') &= 2b^{il}\p^{jk} (2\p(h^2)_{kl}\nabla^2_{ij}\p + \p(h^2)_{ij}\p(h^2)_{kl}) \\
&\quad + 2b^{il} \p^{jk} \left(\gamma_{ij} \gamma_{kl} - 2\a_{ij} \gamma_{kl} + \nabla^2_{ij}\p\nabla^2_{kl}\p\right) \\
&= 2b^{il}\p^{jk} \big(\nabla^2_{ij}\p\nabla^2_{kl}\p +2\p(h^2)_{kl}\nabla^2_{ij}\p + \p(h^2)_{ij}\p(h^2)_{kl} \\
    &\qquad\qquad - 2 \a_{ij} \gamma_{kl} + \gamma_{ij} \gamma_{kl} \big) \\
&= 2b^{il}\p^{jk} \left(\a_{ij}\a_{kl} - 2 \a_{ij} \gamma_{kl} + \gamma_{ij} \gamma_{kl} \right) \\
    &= 2b^{il}\p^{jk} \eta_{ij} \eta_{kl}, \\
(E) - (E') &= c(R_{\beta} - R_{\theta}) \\
    &= cR.\\
\end{align*}
Multiplying everything by \(t\) and adding the result to \eqref{eq:deltchi1} yields the claim.
\end{proof}
We need two more lemmas to obtain a Harnack inequality for $H^{p}$-flow with $0<p\leq 1.$  We start by rewriting the term $R$ in the evolution of \(\chi_2\) when the speed is a function of the mean curvature.

In the sequel, we denote
$F'=\frac{dF}{dH}$
and similarly for higher derivatives.
\begin{lemma}\label{RSphere}
Suppose that $\p=\p(H).$ Then the term $R$ in the evolution equation of $\chi_2$ takes the form
\eq{\label{RSphere1}
R&=2n\fr{\p''\p}{\p'}\br{\Box\p + \p\p^{ij}(h^2)_{ij} - b^{ij}\nabla_i \p \nabla_j \p}-n\fr{\p''\p^2}{\p'}\p^{ij}(h^2)_{ij}\\
		&\hp{=}+2\p^2\p'H+n\br{2\fr{\p''}{\p'}-\fr{\p''^{2}\p}{\p'^{3}}+\fr{\p'''\p}{\p'^{2}}}\p^{ij}\nabla_{i}\p\nabla_j\p\\
 &\hp{=}+\br{\p'H+\p}b^{ir}b^{j}_{r}\nabla_i\p\nabla_j\p-2\p'b^{ij}\nabla_i\p\nabla_j\p.
}
\end{lemma}
\pf{
We calculate the crucial terms in \eqref{R}:
\eq{\Box\operatorname{tr}(\dot{\p})&=\p^{kl}\nabla_{kl}^{2}(\p^{ij}g_{ij})=n\p^{kl}\nabla_{kl}^{2}\p'\\
			&=n\p^{kl}\p'''\nabla_{k}H\nabla_{l}H+n\p^{kl}\p''\nabla_{kl}^{2}H\\
			&=n\frac{\p'''}{\p'^{2}}\p^{kl}\nabla_{k}\p\nabla_{l}\p+n\fr{\p''}{\p'}\Box\p-n\fr{\p''^{2}}{\p'^{3}}\p^{kl}\nabla_{k}\p\nabla_{l}\p\\
			&=n\fr{\p''}{\p'}\Box\p+n\br{\fr{\p'''}{\p'^{2}}-\fr{\p''^{2}}{\p'^{3}}}\p^{kl}\nabla_{k}\p\nabla_{l}\p.}
Furthermore
\eq{\p^{ij,kl}g_{kl}\nabla_{ij}^{2}\p&=n\p''g^{ij}\nabla_{ij}^{2}\p=n\fr{\p''}{\p'}\Box\p}
and
\eq{2F^{kl}\nabla_k\mathrm{tr}(\dot{F})\nabla_lF=2F^{kl}g_{ij}F^{ij,rs}\nabla_k h_{rs}\nabla_lF=2n\frac{F''}{F'}F^{kl}\nabla_kF\nabla_lF.}

Thus
\eq{R&=2n\fr{\p''\p}{\p'}\Box\p+n\br{\fr{\p'''\p}{\p'^{2}}-\fr{\p''^{2}\p}{\p'^{3}}+2\fr{\p''}{\p'}}\p^{kl}\nabla_{k}\p\nabla_{l}\p\\
		&\hp{=}+n\fr{\p''\p^{2}}{\p'}\p^{ij}(h^{2})_{ij}-2n\fr{\p''\p}{\p'}b^{ij}\nabla_{i}\p\nabla_{j}\p+2\p^{2}\p'H\\
		&\hp{=}+(\p^{kl}h_{kl} + \p)b^{ir}b^{j}_{r}\nabla_i \p\nabla_j\p - 2 b^{j}_{k}\p^{kl}\nabla_l\p\nabla_j\p}
and a little rearrangement gives the result.

}
To obtain a Harnack estimate for $H^{p}$-flow, we will have to handle the middle term in \eqref{RSphere1}; this term does not always have the favourable positive sign. To this aim, it is useful to add an auxiliary function of the speed. Using Proposition \ref{thm:Evchi} and \cref{RSphere} it is straightforward to obtain the following evolution equation for $\chi_3=\chi_2+ct\zeta$, where $\zeta=\zeta(\p)$ is a function of $\p.$
\begin{lemma}\label{lemma : lem9}
Let $F=F(H)$. Then under flow \eqref{eq:CurvFlow} the quantity
$\chi_3=\chi_2+ct\zeta$
evolves according to
\eq{\label{Evchibar1}
\del_t \chi_3 -\Box\chi_3 &= \left(\fr{\b-\t}{\d\p} + \p^{ij}(h^2)_{ij} + c\operatorname{tr}(\dot{\p})\right)\chi_2+c\zeta \\
&\hp{=} + t\left(\p^{ij,kl} + 2b^{il}\p^{jk} - \fr{\p^{ij}\p^{kl}}{\d\p}\right)\eta_{ij}\eta_{kl} \\
&\hp{=} + ct\Big\{2n\fr{\p''\p}{\p'}\br{\Box\p + \p\p^{ij}(h^2)_{ij} - b^{ij}\nabla_i \p \nabla_j \p}\\
&\hp{=+t}+\br{\zeta'-n\fr{\p''\p}{\p'}}\p^{ij}(h^2)_{ij}\p+ c\zeta'\operatorname{tr}(\dot{\p})\p+ 2\p^2\p'H\\
&\hp{=+t}+\br{n\br{2\fr{\p''}{\p'}-\fr{\p''^{2}\p}{\p'^{3}}+\fr{\p'''\p}{\p'^{2}}}-\zeta''}\p^{ij}\nabla_{i}\p\nabla_j\p \\
&\hp{=+t} + \br{\p'H+\p}b^{ir}b^{j}_{r}\nabla_i\p\nabla_j\p-2\p'b^{ij}\nabla_i\p\nabla_j\p\Big\}.
}
\end{lemma}

\begin{proof}
\[\partial_t\chi_3-\Box\chi_3=(\partial_t-\Box)\chi_2+c\zeta+ct(\partial_t-\Box)\zeta.\]
Adding \eqref{thm:Evchi1} to
\eq{c\zeta+ct(\partial_t-\Box)\zeta&=c\zeta+ct\left(\zeta'(\partial_t-\Box)F-\zeta''F^{kl}\nabla_{k}F\nabla_lF\right)\\
			&=c\zeta+ct\left(\zeta'FF^{ij}(h^2)_{ij}+c\zeta'FF^{ij}g_{ij}-\zeta''F^{kl}\nabla_kF\nabla_lF\right)}
            gives the result.
\end{proof}

\Cref{thm:Evchi}, and Lemmas \ref{RSphere}, \ref{lemma : lem9} enable us to get a strong Harnack estimate for $H^{p}$-flows; see \Cref{sec:Harnack} and \Cref{thm: main 1}. Due to the presence of $\Box\operatorname{tr}(\dot{\p})$ in $R$ given in \Cref{thm:Evchi}, it is not clear to us whether $\chi_2$ would result in Harnack estimates for curvature flows other than $H^{p}$-flows. As it will be shown, by weakening $\chi_2$ to $\chi_1=\chi_2+tc\p\operatorname{tr}(\dot{\p})$,  we can obtain (weak) Harnack estimates for $p$-powers of 1-homogeneous convex speeds, $0<p\leq 1.$

\begin{proposition}\label{WeakHarnackEv}
For a general curvature function $F$ under flow \eqref{eq:CurvFlow} the quantity
$\chi_1=t(\del_t\p-\t)+\d\p$
satisfies the evolution equation
\eq{\label{WeakHarnackEv12}\del_t\chi_1-\Box\chi_1 &= \br{\frac{\b-\t}{\delta\p} + \p^{ij}(h^2)_{ij} + c\frac{\d-1}{\d} \operatorname{tr}(\dot{\p})}\chi_1\\
&\quad +\fr{c\operatorname{tr}(\dot{\p})\p}{\d}(tc\operatorname{tr}(\dot{\p})+2\d) \\
&\quad + t\p^{ij,kl}\br{\eta_{ij}+c\p g_{ij}}\br{\eta_{kl}+c\p g_{kl}}\\
 &\quad +t\br{2b^{il}\p^{jk}-\fr{\p^{ij}\p^{kl}}{\d\p}}\eta_{ij}\eta_{kl} \\
&\quad + tc\left\{2\p^2\p^{ij}h_{ij}+\left(\br{\p^{ij}h_{ij}+\p}b^{ir}-2\p^{ir}\right)b^{j}_{r}\nabla_{i}\p\nabla_{j}\p\right\}.
}
\end{proposition}
\pf{
We simply use the evolution of $\chi_2,$ cf.~\eqref{thm:Evchi1}, and add it to the evolution of $tc\operatorname{tr}(\dot{\p})\p$. We have
\eq{\label{WeakHarnackEv2}\br{\del_t-\Box}\br{tc\operatorname{tr}(\dot{\p})\p}&=tc\Big(\operatorname{tr}(\dot{\p})\p\p^{ij}(h^2)_{ij}+c\operatorname{tr}(\dot{\p})^2\p+\p\del_t\operatorname{tr}(\dot{\p})\\
                    &\hp{-tk\Big(}-\p\Box\operatorname{tr}(\dot{\p})-2\p^{ij}\nabla_{i}\p\nabla_{j}\operatorname{tr}(\dot{\p})\Big)+c\operatorname{tr}(\dot{\p})\p.
                            }
Note that by \eqref{eq:deltBox} there holds
\eq{\del_t\operatorname{tr}(\dot{\p})=\del_t\br{\p^{ij}g_{ij}}&=\left(\p^{ij,kl}\br{\a_{kl}+c\p g_{kl}}+2FF^{jk}h^i_k\right)g_{ij}-2F^{ij}Fh_{ij}\\
		&=\p^{ij,kl}\br{\a_{kl}+c\p g_{kl}}g_{ij}.}
		
Hence
\eq{\label{WeakHarnackEv3}\del_{t}\chi_{1}-\Box\chi_{1}&= \left(\frac{\b-\t}{\delta\p} +\p^{ij}(h^2)_{ij} + c\operatorname{tr}(\dot{\p})\right)\chi_2 \\
& \quad + t\left(\p^{ij,kl} + 2b^{il}\p^{jk} - \frac{\p^{ij}\p^{kl}}{\delta\p}\right)\eta_{ij}\eta_{kl} + tc R\\
			&\quad+tc\Big(\operatorname{tr}(\dot{\p})\p\p^{ij}(h^2)_{ij}+c\operatorname{tr}(\dot{\p})^2\p+\p\p^{ij,kl}\a_{kl}g_{ij}\\
                    &\hp{-tk\Big(}+c\p^{2}\p^{ij,kl} g_{kl}g_{ij}-\p\Box\operatorname{tr}(\dot{\p})-2\p^{ij}\nabla_{i}\p\nabla_{j}\operatorname{tr}(\dot{\p})\Big)\\
                    &\quad+c\operatorname{tr}(\dot{\p})\p\\
                    		&=\left(\frac{\b-\t}{\delta\p} +\p^{ij}(h^2)_{ij} + c\operatorname{tr}(\dot{\p})\right)\chi_2 \\
& \quad + t\left(\p^{ij,kl} + 2b^{il}\p^{jk} - \frac{\p^{ij}\p^{kl}}{\delta\p}\right)\eta_{ij}\eta_{kl} \\
			&\quad+tc\Big(\operatorname{tr}(\dot{\p})\p\p^{ij}(h^2)_{ij}+c\operatorname{tr}(\dot{\p})^2\p+2\p\p^{ij,kl}\eta_{kl}g_{ij}\\
                    &\hp{-tk\Big(}+c\p^{2}\p^{ij,kl} g_{kl}g_{ij} + 2\p^{2}\p^{ij}h_{ij}\\
&\hp{-tk\Big(}  +(\p^{kl}h_{kl} + \p)b^{ir}b^{j}_{r}\nabla_i \p\nabla_j\p - 2 b^{j}_{k}\p^{kl}\nabla_l\p\nabla_j\p.
\Big)\\
                    &\quad+c\operatorname{tr}(\dot{\p})\p,
}
where we have combined the terms attached to the factor $tc$ and used \eqref{R} and $\eta_{ij}=\a_{ij}-\g_{ij}$. Thus, collecting all terms containing $\eta_{ij}$ and $F^{ij,kl}$ in \eqref{WeakHarnackEv3} we get after some rearranging, that

\eq{\del_{t}\chi_{1}-\Box\chi_{1}&= \left(\frac{\b-\t}{\delta\p} +\p^{ij}(h^2)_{ij} + c\operatorname{tr}(\dot{\p})\right)(\chi_1-tcF\mathrm{tr}(\dot{F})) \\
				&\quad +tc\operatorname{tr}(\dot{\p})\p\p^{ij}(h^2)_{ij}+tc^{2}\operatorname{tr}(\dot{\p})^2\p+c\mathrm{tr}(\dot{\p})\p\\
				&\quad + t\p^{ij,kl}\br{\eta_{ij}+c\p g_{ij}}\br{\eta_{kl}+c\p g_{kl}}\\ 
				&\quad +t\br{2b^{il}\p^{jk}-\fr{\p^{ij}\p^{kl}}{\d\p}}\eta_{ij}\eta_{kl}\\
				&\quad + tc\left(2\p^2\p^{ij}h_{ij}+\left(\br{\p^{ij}h_{ij}+\p}b^{ir}-2\p^{ir}\right)b^{j}_{r}\nabla_{i}\p\nabla_{j}\p\right).
}

It remains to rearrange the first two lines of the previous equation. We have

\eq{&\left(\frac{\b-\t}{\delta\p} +\p^{ij}(h^2)_{ij} + c\operatorname{tr}(\dot{\p})\right)(\chi_1-tcF\mathrm{tr}(\dot{F}))\\
	 &\qquad+ tc\operatorname{tr}(\dot{\p})\p\p^{ij}(h^2)_{ij}+tc^{2}\operatorname{tr}(\dot{\p})^2\p+c\mathrm{tr}(\dot{\p})\p\\
	 =&\left(\frac{\b-\t}{\delta\p} +\p^{ij}(h^2)_{ij} + c\operatorname{tr}(\dot{\p})\right)\chi_1-\left(\frac{\b-\t}{\delta\p} +\p^{ij}(h^2)_{ij}\right)tcF\mathrm{tr}(\dot{F})\\
	 &\qquad + tc\operatorname{tr}(\dot{\p})\p\p^{ij}(h^2)_{ij}+c\mathrm{tr}(\dot{\p})\p\\
	 =&\left(\frac{\b-\t}{\delta\p} +\p^{ij}(h^2)_{ij} + c\operatorname{tr}(\dot{\p})\right)\chi_1-\left(\frac{\chi_{2}}{\delta\p}-1 +t\p^{ij}(h^2)_{ij}\right)cF\mathrm{tr}(\dot{F})\\
	 &\qquad + tc\operatorname{tr}(\dot{\p})\p\p^{ij}(h^2)_{ij}+c\mathrm{tr}(\dot{\p})\p\\
	 =&\left(\frac{\b-\t}{\delta\p} +\p^{ij}(h^2)_{ij} + c\fr{\d-1}{\d}\operatorname{tr}(\dot{\p})\right)\chi_1+\fr{tc^{2}}{\d}F\mathrm{tr}(\dot{F})^{2}+2c\p\mathrm{tr}(\dot{\p}),}
which is the proposed form of the first two lines in \eqref{WeakHarnackEv12}.

}

\section{Harnack Inequalities}
\label{sec:Harnack}

In Euclidean space we recover differential Harnack inequalities for various speeds already discussed in \cite[Corollary 5.11 (1)]{Andrews:09/1994}, as can be seen by evaluating the evolution equation \eqref{thm:Evchi1} with $c=0.$

\begin{remark}\label{HarnackEucCor}
 Let $x$ be a strictly convex solution of \eqref{eq:CurvFlow} in Euclidean space, i.e. $c=0$. Note that the evolution of $\chi_{1}$ simplifies tremendously now:
 \eq{
\partial_t \chi_1 -\Box\chi_1 &= \left(\frac{\b-\t}{\delta\p} +\p^{ij}(h^2)_{ij}\right)\chi_1 \\
& \quad + t\left(\p^{ij,kl} + 2b^{il}\p^{jk} - \frac{\p^{ij}\p^{kl}}{\delta\p}\right)\eta_{ij}\eta_{kl}.
}

For curvature functions $F$ and suitable $\delta$ so that the second term is non-negative, we obtain a Harnack inequality. For example, let $f$ be a $1$-homogeneous, inverse concave curvature function; that is, the curvature function
 \eq{\label{inverse}\~f(\mc{W}):=\fr{1}{f(\mc{W}^{-1})}}
 is concave. In \cite[p.~112]{Urbas:/1991} it is shown that in this case $f$ satisfies
 \eq{(f^{ij,kl}+2f^{ik}b^{jl})\eta_{ij}\eta_{kl}\geq 2f^{-1}f^{ij}f^{kl}\eta_{ij}\eta_{kl}}
 for all symmetric matrices $(\eta_{ij})$. For $\a\neq 0,$ setting
 \eq{F=\fr{|\a|}{\a}f^{\a},}
 we obtain
 \eq{F^{ij}=|\a|f^{\a-1}f^{ij},\quad F^{ij,kl}=(\a-1)|\a|f^{\a-2}f^{ij}f^{kl}+|\a|f^{\a-1}f^{ij,kl}}
 and hence, in the sense of bilinear forms,
 \eq{F^{ij,kl}+2b^{il}F^{jk}-\fr{F^{ij}F^{kl}}{\d F}&=|\a|f^{\a-1}\Big(f^{ij,kl}+\fr{\a-1}{f}f^{ij}f^{kl}\\
 			&\qquad\qquad+2b^{il}f^{jk}-\fr{\a}{\d f}f^{ij}f^{kl}\Big)\\
 				&\geq |\a|f^{\a-2}\br{\a+1-\fr{\a}{\d}}f^{ij}f^{kl}\\
				&\geq 0,}
if we chose $\d\geq\fr{\a}{\a+1}$ when $\a \geq 0$ and $\d \leq \fr{\a}{\a+1}$ when $\a \leq 0.$ However, in order to apply this estimate in a maximum principle argument for $\chi_{1},$ we need that $\chi_{1}$ to be positive initially and hence $F$ and $\d$ must have the same sign. Thus the only allowed pairs $(F,\d)$ are
 
$$\p(f)=f^{\a},\quad 0<\a<\8, \quad \d \geq \fr{\a}{\a+1},$$
for the contracting flows and for expanding flows we can allow the cases
$$\p(f)=-f^{-\b}, \quad 0<\b<1,\quad \d \leq \fr{\b}{\b-1}.$$
In those situations we obtain
 $$\fr{\chi_{1}}{t}=\partial_t \p-b^{ij}\nabla_i\p\nabla_j\p+\frac{\d\p}{t}>0.$$
Compare with \cite[Theorem 5.6, Corollary 5.11]{Andrews:09/1994}. To the best of our knowledge, this is the first time these Harnack inequalities (even in Euclidean space) have been proved in such generality in the parametric setting.
\end{remark}
Now we move on to the spherical case. Before we can prove the main theorem, for convenience we provide the proof of an inequality for curvature functions, the idea of which can be found in \cite[Theorem 2.3]{Andrews:/2007}.

\begin{lemma}\label{f-Lemma}
Let $f=f(h_{ij},g_{ij})$ be a monotone, $1$-homogeneous curvature function defined on $\Gamma_+.$
Then 
\begin{equation}\label{1}\left(f^{ik}b^{jl}-\frac{f^{ij}f^{kl}}{f}\right)\eta_{ij}\eta_{kl}\geq 0
\end{equation}
for all symmetric matrices $\eta_{ij}.$
\end{lemma}

\begin{proof}
Take a coordinate system, such that
$$h_{ij}=\kappa_i \delta_{ij}.$$
Then we also have
$$f^{ij}=\frac{\partial f}{\partial h_{ij}}=\frac{\partial f}{\partial \kappa_i}\delta^{ij}\equiv f^i\delta^{ij}$$
Thus the left hand side of \eqref{1} becomes
$$\frac{f^i}{\kappa_j}\eta_{ij}^2-\frac{f^i\eta_{ii}f^j\eta_{jj}}{f}\geq \frac{f^i}{\kappa_i}\eta^2_{ii}-\frac{f^i\eta_{ii}f^j\eta_{jj}}{f},
$$
where we have just thrown away all non-diagonal entries of $\eta_{ij}$ (Note that this is not a waste, since we have to prove the thing anyway for all matrices.).
From now on we denote $(\eta_{ii})$ simply by $\eta_i.$ Thus we have to prove
$$\forall \eta\in\mathbb{R}^n\colon\frac{f^i}{\kappa_i}\eta^2_{i}-\frac{f^i\eta_{i}f^j\eta_{j}}{f}\geq 0
$$
Define the $(n-1)$ dimensional linear subspace
$$S=\{(\xi_i)\in \mathbb{R}\colon f^i\xi_i=0 \}.
$$
Since $f^i\kappa_i=f=f(\kappa)>0,$ we have
$$\mathbb{R}^n=S\oplus \langle \kappa\rangle.$$
Thus
$$\eta=\xi+a\kappa, \quad \xi\in S,\quad a\in \mathbb{R},$$
and we may assume $a=1,$
for if $a=0$ there is nothing to prove and if $a\neq 0$ take $\tilde{\eta}=\frac{\eta}{a}.$
The desired inequality becomes, due to the homogeneity,
$$\frac{f^i}{\kappa_i}(\xi_i+\kappa_i)^2-f=\frac{f^i}{\kappa_i}\xi_i^2+2\frac{f^i}{\kappa_i}\xi_i\kappa_i=\frac{f^i}{\kappa_i}\xi_i^2\geq 0.$$
\end{proof}

\begin{theorem}
Let $f$ be a strictly monotone, $1$-homogeneous, convex curvature function, $0<p\leq 1$, and let $\p=f^p$. Then under flow \eqref{eq:CurvFlow} with $c\ge 0$, $\chi_{1}$ satisfies
$$\frac{\chi_{1}}{t}=\partial_t \p-b^{ij}\nabla_i\p\nabla_j\p+\frac{p\p}{(p+1)t}>0\quad \forall t\in (0,T).$$
\end{theorem}
\begin{proof}
In view of the maximum principle and that $\chi_1$ is manifestly positive at $t=0$, it suffices to show that the right-hand side of (\ref{WeakHarnackEv12}) is positive whenever at some point in space-time $\chi_1=0$. 
Due to \cref{f-Lemma} there holds
\eq{f^{ik}b^{jl}\eta_{ij}\eta_{kl}\geq f^{-1}\br{f^{ij}\eta_{ij}}^2}
for all symmetric matrices $\eta$ and hence
\eq{F^{ik}b^{jl}\eta_{ij}\eta_{kl}=pf^{p-1}f^{ik}b^{jl}\eta_{ij}\eta_{kl}\geq pf^{p-2}(f^{ij}\eta_{ij})^{2}\geq \fr 1p F^{-1}\br{F^{ij}\eta_{ij}}^2.}
Hence for $\d=\fr{p}{p+1}$ we have
\eq{\label{convex F}\br{2b^{il}F^{jk}-\fr{F^{ij}F^{kl}}{\d F}}\eta_{ij}\eta_{kl}&\geq \fr{1-p}{p}F^{-1}F^{ij}F^{kl}\eta_{ij}\eta_{kl}\\
			&=\fr{1-p}{p}F^{-1}F^{ij}F^{kl}(\eta_{ij}+cFg_{ij})(\eta_{kl}+cFg_{kl})\\
            &\hp{=}-2c\fr{1-p}{p}F^{ij}\eta_{ij}F^{kl}g_{kl}-c^2\fr{1-p}{p}F\br{F^{ij}g_{ij}}^2\\
            &=\fr{1-p}{p}F^{-1}F^{ij}F^{kl}(\eta_{ij}+cFg_{ij})(\eta_{kl}+cFg_{kl})\\
            &\hp{=}-\fr{2c}{t}\fr{1-p}{p}\chi_{1}F^{kl}g_{kl}+\fr{2\d c}{t}\fr{1-p}{p}FF^{kl}g_{kl}\\
            &\hp{=}+c^2\fr{1-p}{p}F\br{F^{ij}g_{ij}}^2,
            }
where in the last equality we have used
\eq{F^{ij}\eta_{ij}=\b-\t=\fr{\chi_{1}-\d F}{t}-cFF^{ij}g_{ij}.}
The first term in the last equality of \eqref{convex F} when added to the term involving $F^{ij,kl}$ in \eqref{WeakHarnackEv12} produces a positive term:
\eq{F^{ij,kl}+\fr{1-p}{p}F^{-1}F^{ij}F^{kl}&=pf^{p-1}f^{ij,kl}\geq 0}
as bilinear forms due to the convexity of $f$.
 The other terms in \eqref{convex F} do no harm in applying the maximum principle.
On the other hand, note that any strictly monotone, 1-homogeneous curvature function $f$ satisfies $f b^{ij}\geq f^{ij}.$ Therefore
\eq{\left(\br{\p^{ij}h_{ij}+\p}b^{ir}-2\p^{ir}\right)b^{j}_{r}\nabla_{i}\p\nabla_{j}\p&=f^{p-1}\br{(p+1)-2p}f^{ir}b^{j}_{r}\nabla_{i}F\nabla_{j}\geq  0.} 
\end{proof}

Employing the evolution \cref{Evchibar1}, we can obtain a stronger Harnack inequality for the speed \(\p = H^{p}\) with \(p \in (0,1)\); case $p=1$ was considered in \cite{BryanIvaki:08/2015}.

\begin{theorem} \label{thm: main 1}
Consider a solution of \eqref{eq:CurvFlow} with \(\p = H^{p}\) and $c\geq 0.$ If $\frac{1}{2}+\frac{1}{2n}\leq {p}< 1,$ then
\[
\partial_t H^{p} - b^{ij}\nabla_iH^{p}\nabla_jH^{p} - \frac{c {p}}{2{p}-1}H^{2{p}-1} + \frac{{p}}{{p}+1} \frac{H^{p}}{t} > 0.
\]
If $0<{p}\leq \frac{1}{2} + \frac{1}{2n}$ or $p=1$, then

\[
\partial_t H^{p} - b^{ij}\nabla_iH^{p}\nabla_jH^{p} - c n{p}H^{2{p}-1} + \frac{{p}}{{p}+1} \frac{H^{p}}{t} > 0.
\]
\end{theorem}
\begin{proof}
In order to prove Theorem \ref{thm: main 1}, we need to show that for
$$\p=H^{p},\quad \d=\fr{p}{p+1},$$
the quantity $\chi_3$ preserves its positivity at all $t>0.$ Here $\zeta$ is chosen to be
$$\zeta(\p)=\begin{cases} p\br{n-\fr{1}{2p-1}}\p^{2-\fr{1}{p}}, &\fr{1}{2}+\fr{1}{2n}<p< 1\\
               0, &0<p\leq \fr{1}{2}+\fr{1}{2n}~\text{or}~p=1.\end{cases}$$
However, to avoid confusion, we will keep the general form as long as possible.
At time $t=0,$ $\chi_3$ is positive.
Thus suppose there exists a first time $t_0$ and a point $x_0$ in $M_{t_0},$ such that $\chi_3(t_0,x_0)=0.$
Then we also obtain
$$\chi_2=-ct\zeta\Rightarrow \b-\t=-\fr{\d\p}{t}-c\zeta.$$
Thus, using \eqref{Evchibar1} and $\Box\p+\p\p^{ij}(h^{2})_{ij}-b^{ij}\nabla_{i}\p\nabla_{j}\p=\b-\t$, we obtain at $(t_0,x_0)\colon$
\eq{&0\geq \del_t\chi_3-\Box\chi_3\\
        &=2c\zeta-2cn\d\fr{\p''\p^2}{\p'}+t\br{\p^{ij,kl}+2b^{il}\p^{jk}-\fr{\p^{ij}\p^{kl}}{\d\p}}\eta_{ij}\eta_{kl}\\
        &\hp{=}+ct\Big\{\fr{c\zeta^2}{\d\p}-2cn\fr{\p''\p}{\p'}\zeta+2\p^2\p'H+c\br{\zeta'\p-\zeta}\p^{ij}g_{ij}    \\
        &\hp{=+t} +\br{\zeta'\p-n\fr{\p''\p^2}{\p'}-\zeta}\p^{ij}(h^2)_{ij}+\br{\p'H+\p}b^{ir}b^j_r\nabla_i\p\nabla_j\p\\
        &\hp{=+t}+\br{n\br{2\fr{\p''}{\p'}-\fr{\p''^{2}\p}{\p'^{3}}+\fr{\p'''\p}{\p'^{2}}}-\zeta''}\p^{ij}\nabla_{i}\p\nabla_j\p\\
        &\hp{=+t}-2\p'b^{ij}\nabla_i\p\nabla_j\p\Big\}\\
        &\geq2c\zeta-2cn\d\fr{\p''\p^2}{\p'}+t\br{\p^{ij,kl}+2b^{il}\p^{jk}-\fr{\p^{ij}\p^{kl}}{\d\p}}\eta_{ij}\eta_{kl}\\
        &\hp{=}+ct\Big\{\fr{c\zeta^2}{\d\p}+2\p^2\p'H-2cn\fr{\p''\p}{\p'}\zeta+cn\br{\zeta'\p-\zeta}\p'    \\
        &\hp{=+t} +\br{\zeta'\p-n\fr{\p''\p^2}{\p'}-\zeta}\p^{ij}(h^2)_{ij}\\
        &\hp{=+t}+\br{n\br{2\fr{\p''}{\p'}-\fr{\p''^{2}\p}{\p'^{3}}+\fr{\p'''\p}{\p'^{2}}}-\zeta''+\fr{\p}{\p'H^2}-\fr{1}{H}}\p^{ij}\nabla_{i}\p\nabla_j\p\Big\},}
        
where we used the estimate

\eq{&(\p'H+\p)b^{ir}b^{j}_{r}\nabla_{i}F\nabla_{j}\p-2\p'b^{ij}\nabla_{i}F\nabla_{j}F\\
	\geq &(p+1)\fr{F}{H}b^{ij}\nabla_{i}F\nabla_{j}F-2p\fr{F}{H}b^{ij}\nabla_{i}F\nabla_{j}F\\
	=&(1-p)\fr{F}{H}b^{ij}\nabla_{i}F\nabla_{j}F\\
	\geq &(1-p)\fr{F}{F'H^{2}}F^{ij}\nabla_{i}F\nabla_{j}F. }

To finish the proof, we need to show that the right-hand side is positive. If $\zeta=0$, this is straightforward:
\[\p''<0,\quad n\br{2\fr{\p''}{\p'}-\fr{\p''^{2}\p}{\p'^{3}}+\fr{\p'''\p}{\p'^{2}}}+\fr{\p}{\p'H^2}-\fr{1}{H}\geq 0\]
and
$$\br{\p^{ij,kl}+2b^{il}\p^{jk}-\fr{\p^{ij}\p^{kl}}{\d\p}}\eta_{ij}\eta_{kl}\geq0\quad\mbox{for}\quad\delta=\frac{p}{p+1}.$$
For the second case that $\zeta\neq 0$, note that
$$2c\zeta-2cn\d\fr{\p''\p^2}{\p'}\geq 0\quad\mbox{for}\quad p\geq \frac{n+1}{2n},$$
$$\br{\p^{ij,kl}+2b^{il}\p^{jk}-\fr{\p^{ij}\p^{kl}}{\d\p}}\eta_{ij}\eta_{kl}\geq0\quad\mbox{for}\quad\delta=\frac{p}{p+1},$$
$$\fr{c\zeta^2}{\d\p}-2cn\fr{\p''\p}{\p'}\zeta+cn\br{\zeta'\p-\zeta}\p'\geq 0\quad\mbox{for}\quad p\geq\frac{n+1}{2n},~\delta=\frac{p}{p+1},$$
$$\zeta'\p-n\fr{\p''\p^2}{\p'}-\zeta\geq 0\quad\mbox{for}\quad\frac{n+1}{2n}\leq p\leq 1,$$
$$n\br{2\fr{\p''}{\p'}-\fr{\p''^{2}\p}{\p'^{3}}+\fr{\p'''\p}{\p'^{2}}}-\zeta''+\fr{\p}{\p'H^2}-\fr{1}{H}=0.$$
\end{proof}

\section{Preserving convexity}
\label{sec:convexity}

In the derivation of the Harnack inequalities we have assumed the strict convexity of flow hypersurfaces. In this section, we show that this assumption is justified by proving strict convexity is preserved for all flows in the sphere for which we could prove the Harnack inequality. In Euclidean space, the question of preserved convexity has been addressed more thoroughly. It is also known that there is a variety of examples where convexity is lost for contracting flows \cite{AndrewsMcCoyZheng:07/2013}; the authors also discuss necessary and sufficient conditions for preserving convexity. In other special situations preserved convexity was proved, e.g., see \cite{Andrews:/1994b,Andrews:/2007,Andrews:04/2010,Schulze:/2005}.
\begin{proposition}
 Let $M_{0}\sub\S^{n+1}$ be a closed and strictly convex hypersurface. Suppose that $f\in C^{\8}(\G_{+})\cap C^{0}(\bar{\G}_{+})$ is a strictly monotone, $1$-homogeneous and convex curvature function and let $0<p<\8$. Let $x$ be the solution to \eqref{eq:CurvFlow} with $F=f^{p}$ and with initial hypersurface $M_{0}$. Then all flow hypersurfaces $M_{t}=x(M,t)$ are strictly convex.
\end{proposition}
\pf{

Let $T$ be the first time, where the strict convexity is lost. Then on the time interval $[0,T)$ the dual flow defined via the Gauss map is well defined and reads
\eq{\label{DualFlow}\dot{\~x}=\fr{1}{\~f^{p}}\~\nu,}
where $\~f$ is the inverse curvature function defined in \eqref{inverse} and $\~f$ is now evaluated at $\~\k_i=\k_i^{-1};$ see \cite{Gerhardt:/2015} for the derivation of the dual flow. Due to the properties of $f$, $\~f$ is 1-homogeneous, strictly monotone, concave and vanishes on the boundary of $\G_{+} $; see \cite[Lemma 2.2.12, Lemma 2.2.14]{Gerhardt:/2006}. For flows of the kind \eqref{DualFlow} uniform curvature estimates were deduced in \cite[Lemma 4.7]{MakowskiScheuer:/2013}, implying that the $\~\k_{i}$ are bounded. This means that up to time $T$ uniform convexity is preserved for the original flow, which contradicts the definition of $T,$ if $T$ is not the collapsing time.
}

\bibliographystyle{amsplain}
\bibliography{Bibliography.bib}
\end{document}